\documentclass[12 pt,reqno]{amsart}

\usepackage{latexsym, amsmath, amssymb, longtable, booktabs,amscd,microtype,booktabs,cases}
\usepackage[square,numbers]{natbib}
\usepackage{graphicx}
\usepackage[dvipsnames]{xcolor}
\usepackage{caption}
\usepackage{dsfont}
\usepackage[all]{xy}
\usepackage{enumerate}
\usepackage{upgreek}
\usepackage{bm}
\usepackage{relsize}
\usepackage{multirow}
\usepackage{tikz}
\usetikzlibrary{matrix,arrows,decorations.pathmorphing}
\usepackage{collectbox}
\usepackage{enumerate}
\usepackage[colorinlistoftodos,prependcaption, textsize=tiny]{todonotes}
\reversemarginpar
\usepackage{amsmath}
\usepackage{enumitem}
\usepackage{hyperref}
\usepackage{cleveref}
\let\oldtocsection=\tocsection

\let\oldtocsubsection=\tocsubsection

\let\oldtocsubsubsection=\tocsubsubsection

\renewcommand{\tocsection}[2]{\hspace{0em}{\vspace{0.5em}}\oldtocsection{#1}{#2}}
\renewcommand{\tocsubsection}[2]{\hspace{1em}{\vspace{0.5em}}\oldtocsubsection{#1}{#2}}
\renewcommand{\tocsubsubsection}[2]{\hspace{2em}\oldtocsubsubsection{#1}{#2}}
\hypersetup{
	colorlinks,
	citecolor=red,
	filecolor=black,
	linkcolor=blue,
	urlcolor=black
}

\numberwithin{equation}{section}
\newcommand{\Z}{\mathbb{Z}}

\newcommand{\al}{\alpha}

\newcommand{\lm}{\underbar{\textbf{m}}}

\newcommand{\C}{\mathbb{C}}

\newcommand\N{\mathbb{N}}

\newcommand\U{\mathcal{U}}

\def\gg{{\mathfrak{g}}}

\newcommand\M{{\mathcal{M}}}

\newcommand{\m}{\textbf{m}}

\newcommand{\lr}{\underbar{\textbf{r}}}

\makeatother

\newtheorem{thm}{Theorem}[section]
\newtheorem{theorem}[thm]{Theorem}

\newtheorem{prop}[thm]{Proposition}

\theoremstyle{definition}

\newtheorem{remark}[thm]{Remark}

\theoremstyle{definition}

\theoremstyle{remark}

\theoremstyle{remark}

\makeatletter
\def\imod#1{\allowbreak\mkern10mu({\operator@font mod}\,\,#1)}
\makeatother
\title[Classification of Harish-Chandra modules]{\textbf{Classification of irreducible Harish-Chandra modules over map full toroidal Lie algebras}}
\author{Sudipta Mukherjee}
\begin{document}
	\maketitle
	\begin{abstract}
		A natural higher dimensional analogue of the affine-Virasoro algebra is the full toroidal Lie algebra. In this paper, we classify irreducible Harish-Chandra modules for map full toroidal Lie algebras. We show that every such module is either a cuspidal or a highest weight module. Furthermore, we prove that they turn out to be single point evaluation modules.
		\end{abstract}
	
	\section{Introduction:} The study of infinite-dimensional Lie algebras and their representations has played a central role in various areas of mathematics and mathematical physics, due to their wide-ranging applications in conformal field theory, string theory, and soliton theory. A widely explored class of representations consists of those that decompose into finite-dimensional weight spaces. Such representations are commonly referred to as Harish-Chandra modules in the literature, and their classification for various classes of infinite-dimensional Lie algebras has been the subject of extensive study.\par
	
	 Among the most commonly studied infinite-dimensional Lie algebras is the affine Kac-Moody algebra as it generalizes the theory of finite dimensional simple Lie algebras. Affine Lie algebras can be constructed as the universal central extension of the loop algebra associated with a finite dimensional simple Lie algebra, i.e., the Lie algebra of maps from the unit circle into a finite-dimensional simple Lie algebra. \par
	 
	Another fundamental infinite-dimensional Lie algebra is the Virasoro algebra, which plays an important role in theoretical physics. Virasoro Lie algebras can be realized as the central extension of the Lie algebra of derivations of the Laurent polynomial ring in one variable. The classification of all irreducible modules with finite dimensional weight spaces over the Virasoro algebra was conjectured by V. Kac \cite{bom}, and subsequently proved by O. Mathieu in 1992 \cite{om}.\par
	
	 The Virasoro algebra $\mathrm{Vir}$ acts naturally on the derived algebra of an affine Kac-Moody algebra by derivation. Since they are closely related, it naturally suggests that they can be considered simultaneously, that is, as one algebraic structure. As a result, the so-called affine-Virasoro algebra naturally appears, which is the semidirect product of the Virasoro algebra and an affine Kac-Moody Lie algebra. The representation theory of affine-Virasoro algebras has attracted significant attention due to its connections with conformal field theory, number theory, and soliton theory. Affine-Virasoro algebras with a common center have been studied by many mathematicians \cite{gao,ji, ka,bo}.\par
	
	The higher dimensional analogue of the Virasoro algebra is the Witt algebra. Let $A_n$ denote the Laurent polynomial ring in $n$ variables. The Lie algebra of derivations of \( A_n \), denoted by $\mathcal{W}_n$, is known as the Witt algebra. Unlike the one variable case,  $\mathcal{W}_n$ is centrally closed for $n \geq 2$. The Witt algebra, as a classical and well-studied example of an infinite-dimensional Lie algebra, has attracted substantial attention in the literature~\cite{bill,guo,guo11, rao2}. In \cite{bill}, Billig and Futorny classified all irreducible modules with finite dimensional weight spaces over $\mathcal{W}_n$.\par
	
	Similarly, toroidal Lie algebras generalize affine Lie algebras to higher dimensions and are defined as the universal central extensions of the multiloop algebra of finite dimensional simple Lie algebras. In contrast to the affine case, where the center is one-dimensional, toroidal Lie algebras possess an infinite-dimensional center, making their representation theory significantly more intricate and interesting. Using representations of the toroidal algebras, one can construct hierarchies of non-linear PDEs \cite{yb2}. The classification of irreducible integrable Harish-Chandra modules in this setting was carried out by Eswara Rao \cite{2004}. \par
	
	Full toroidal Lie algebras are the natural multivariable generalizations of affine-Virasoro algebras. They are constructed as the semidirect product of the derived algebra of toroidal Lie algebra and Witt algebra. Due to the nontrivial action of \( \mathcal{W}_n \) on the center \( \mathcal{K} \), the center of the full toroidal Lie algebra becomes finite-dimensional. Very recently, S. Pal provided a classification of all irreducible Harish-Chandra modules over full toroidal Lie algebras \cite{sou} . \par
	
	In recent years, the representation theory of map algebras has emerged as an active and rapidly developing area of research. Let $B$ be a finitely generated commutative associative unital algebra. For a Lie algebra $\mathfrak{g}$, the tensor product $\mathfrak{g}(B) := \mathfrak{g} \otimes B$ is called the map algebra associated with $\mathfrak{g}$. The map algebras corresponding to $B = \mathbb{C}[t, t^{-1}]$ and $B = A_n$ are called the loop algebra and multiloop algebra, respectively. For $B=\C[t,t^{-1}]$, Harish-Chandra modules over $\mathrm{Vir}(B)$ were classified in \cite{ZRX}. This result was later generalized for any finitely generated commutative associative unital algebra $B$ by A. Savage in \cite{SAV}. Very recently, Eswara Rao et al. classified irreducible integrable Harish-Chandra modules for map affine-Virasoro algebras \cite{fi}.  
	\par

	  Let $\uptau$ denote the full toroidal Lie algebra, and let $B$ be a finitely generated commutative associative unital algebra. In this paper, we study the irreducible representations of the  map full toroidal Lie algebra $\uptau(B)$ which has finite dimensional weight spaces. We show that all such modules are either cuspidal or highest weight modules. We also prove that in both cases they are single point evaluation modules (see Section 1 for definition). In other words, these modules are irreducible modules for $\uptau$. 	Our classification relies on a recent result by S. Pal \cite{sou}, who classified irreducible Harish-Chandra modules for $\uptau$. Using this, we establish our main theorem~(\Cref{fc}), completing the classification in the more general setting of map full toroidal Lie algebras.
	  
	\subsection{Organization of the paper.} This paper aims to classify all irreducible Harish-Chandra modules over the map full toroidal Lie algebra.
	 Section 2 lays the foundation for the paper. We begin by recalling the definition and structure of the full toroidal Lie algebra $\uptau$, along with its map version $\uptau(B)$. We also discuss the representation theory of toroidal Lie algebras that form the basis of our classification results.\par
	
	Section 3 is devoted to the classification of irreducible cuspidal Harish-Chandra modules over $\uptau(B)$. We begin by proving that all central operators in $\mathcal{K}(B)$ act trivially on $V$ (Theorem \ref{po}). Subsequently, we establish the existence of a cofinite ideal $J$ of $B$ such that $\uptau(J)V=0$ (Theorem \ref{va}). Propositions \ref{vb} and \ref{vc} further refine this result by showing that the ideal $J$ can, in fact, be taken to be a maximal ideal of $B$. Using these results, and applying the classification theorem from \cite{sou}, we arrive at the main result of this section (Theorem~\ref{ds}), which provides a complete description of irreducible cuspidal Harish-Chandra modules over $\uptau(B)$.\par
	
	In Section 4, we classify irreducible Harish-Chandra modules over $\uptau(B)$ whose weight spaces are not uniformly bounded. We show that all such modules are necessarily of highest weight type (\Cref{oi}). These modules are constructed as irreducible quotients of generalized Verma modules, induced from irreducible cuspidal modules over some suitable subalgebra of $\uptau(B)$. We further show that these modules are also single point evaluation modules (\Cref{sx}), and we provide an explicit description of their highest weight spaces (\Cref{xs}). Section 5 brings together the results from the previous sections to present a complete classification of all irreducible Harish-Chandra modules over $\uptau(B)$ (\Cref{fc}).
	\section{Preliminaries}
	\subsection{Notations:} Throughout this paper, we consider all vector spaces, algebras, and tensor products to be over the field of complex numbers $\C$. We shall denote the set of integers, natural numbers, and real numbers by $\Z, \N, \mathbb{R}$ respectively. We write elements of $\Z^n$ and $\C^n$ in boldface. For an $n$-tuple ${\bf{m}}=(m_1,m_2,\ldots,m_n) \in \Z^n$, we denote $(m_2,\ldots,m_n) \in \Z^{n-1}$ by $\lm$. For two elements ${\bf{m}}, {\bf{k}} \in \Z^n$, we say that ${\bf{m}} \geq {\bf{k}}$ if $m_i \geq k_i$ for all $1 \leq i \leq n$. Let $(-|-)$ denote the standard inner product on $\C^n$. The universal enveloping algebra of a Lie algebra $\gg$ is denoted by $\mathcal{U}(\gg)$. Let $\{{\bf{e_1}}, {\bf{e_2}},\ldots,{\bf{e_n}} \}$ denote the standard basis of $\C^n$. The dual space of a vector space $V$ will be denoted by $V^*$. The notation $A^{\times}$ stands for the set $A \setminus \{0\}$.
	\subsection{Basics of Full Toroidal Lie Algebra}
	
	Let $\gg$ be a finite dimensional Lie algebra with Cartan subalgebra $\mathfrak{h}$. It is well known that $\gg$ admits a non-degenerate symmetric bilinear form $\left\langle \cdot,\cdot \right\rangle$, which is unique up to a scalar.\par
	 For a positive integer $n \geq 3$, let $A_n= \C[t_1^{\pm1}, t_2^{\pm1}, \ldots, t_n^{\pm1}]$ denote the Laurent polynomial ring in $n$ variables $t_1, t_2, \ldots, t_n$. For an $n$-tuple  ${\bf{m}}= (m_1, m_2, \ldots, m_n) \in \Z^{n}$, define  $t^{\bf{m}}=t_1^{m_1} t_2^{m_2} \cdots t_n^{m_n} \in A_n$. The multiloop algebra is defined as $\mathcal{L}(\gg)= \gg \otimes A_n$, with the Lie bracket
	 \begin{center}
	 	$[x \otimes t^{{\bf{m}}}, y \otimes t^{{\bf{k}}}] = [x,y] \otimes t^{{\bf{m+k}}}$, \quad $\forall x,y \in \gg, \,\,\,\, {\bf{m}}, {\bf{k}} \in \Z^{n}$.
	 \end{center}
 Define differential elements $K_i= t_i^{-1} dt_i$\, for $1 \leq i \leq n$ and let $\Omega_{A_n}$ be the free left $A_n$-module generated by $\{K_1, K_2, \ldots, K_n\}$. The differential map $d: A_n \rightarrow \Omega_{A_n}$ is given by $d(t^{\textbf{m}})= \sum_{i=1}^{n} m_i t^{\textbf{m}} K_i$. More precisely, we have
 \medskip
 \begin{center}
 	$\Omega_{A_n}= \mathrm{Span}\, \{t^{\textbf{m}}K_{i}: \textbf{m} \in \Z^n,\, 1\leq i \leq n\},\quad \text{Im}(d)=\text{Span}\, \{ \sum_{i=1}^n m_{i}\,t^{\textbf{m}}K_{i}: \m \in \Z^n\}.$
 \end{center}
 The quotient space $\mathcal{K}=\Omega_{A_n}/\text{Im}\mathrm(d)$ is called the space of  K\"ahler differentials, and is explicitly given by
 \begin{center}
 	$\mathcal{K}= \text{Span}\, \{t^{\textbf{m}}K_{i}: \textbf{m} \in \Z^n,\,\,\, 1\leq i \leq n,\,\,\, \sum_{i=1}^n m_{i}\,t^{\textbf{m}}K_{i}=0\}$.
 \end{center} 
 We use the same notation for the image of $t^{\m} K_i$ in $\mathcal{K}$. It is known (see \cite{kl, yokunama}) that the universal central extension of  $L(\gg)$ is given by $\hat{\mathcal{L}}(\gg)= \mathcal{L}(\gg) \oplus \mathcal{K}$, with the Lie brackets:
 \begin{equation}\label{a1}
 [x \otimes t^{{\bf{m}}},\, y \otimes t^{{\bf{k}}}] = [x,y] \otimes t^{{\bf{m+k}}} + \sum_{i=1}^{n} m_i\, t^{{\bf{m+k}}} K_i
 \end{equation}
\begin{equation}\label{a2}
\hspace{2cm} [x \otimes t^{\bf{m}}, z]= [z, z']=0, \hspace{0.5cm} \forall x,y \in \gg, \quad {\bf{m}}, {\bf{k}} \in \Z^{n},\, z, z' \in \mathcal{K}.
\end{equation}

\par
 The final component in the definition of a full toroidal Lie algebra is the derivation algebra of $A_n$, denoted by $\mathcal{W}_n$. For $1 
 \leq i \leq n$, let $d_i= t_i  \frac{\partial }{\partial t_i}$, which acts on $A_n$ by derivation. Then we know that $\mathcal{W}_n$ is generated by the set $\{t^{\bf{m}} d_i: {\bf{m}} \in \Z^n,\, 1 \leq i \leq n\}$. This infinite dimensional derivation algebra is a very well known classical object and is popularly known as the Witt algebra of rank $n$. Lie bracket on $\mathcal{W}_n$ is given by 
 \begin{center}
 	$[t^{\bf{m}} d_i, t^{\bf{k}} d_j]= k_i\, t^{\bf{m+k}} d_j - m_j\,t^{\bf{m+k}} d_i$. 
 \end{center}
$\mathcal{W}_n$ has a unique extension on $\hat{\mathcal{L}}(\gg)$ by the action
\begin{center}
	$t^{\bf{m}}d_i\, (x \otimes t^{\bf{k}})= k_i\, x \otimes t^{\bf{m+k}}$,\\ \vspace{0.25cm} $t^{\bf{m}}d_i\,(t^{\bf{k}} K_j)= k_i t^{\bf{m+k}} K_j + \delta_{i,j} \sum_{p=1}^{n} m_p t^{\bf{m+k}} K_p$.
\end{center}
Moreover, $\mathcal{W}_n$ admits two non-trivial 2-cocycles with values in $\mathcal{K}$ (see \cite{yb}), given by
\begin{align*}
	\phi_1(t^{\bf{m}} d_i, t^{\bf{k}} d_j)= - k_i m_j \sum_{p=1}^{n} m_p t^{\bf{m+k}} K_p, \quad \quad \phi_2(t^{\bf{m}} d_i, t^{\bf{k}} d_j)= m_i k_j \sum_{p=1}^{n} m_p t^{\bf{m+k}} K_p.
\end{align*}
We are now ready to define the full toroidal Lie algebra. Let $\phi$ be an arbitrary linear combination of $\phi_1$ and $\phi_2$. Then, the full toroidal Lie algebra associated with $\gg$ and $\phi$ is defined as 
\begin{center}
	$\uptau= \mathcal{L}(\gg) \oplus \mathcal{K} \oplus \mathcal{W}_n$.
\end{center}
So, as a vector space
\begin{center}
	$\uptau= \gg \otimes \C[t_1^{\pm1}, t_2^{\pm1}, \ldots, t_n^{\pm1}]  \oplus \sum\limits_{\substack{i=1 \\ \textbf{m} \in \Z^{n}}}^n \mathbb{C} t^{\textbf{m}} K_i \hspace{0.1cm} \bigoplus\limits_{\substack{i=1 \\ \textbf{m} \in \Z^{n}}}^n \mathbb{C} t^{\textbf{m}} d_i$.
\end{center}
The bracket relations in $\uptau$ include \eqref{a1}, \eqref{a2} and
\begin{equation*}
	[t^{\bf{m}}d_i, t^{\bf{k}} K_j]= k_i t^{\bf{m+k}} K_j + \delta_{i,j} \sum_{p=1}^{n} m_p t^{\bf{m+k}} K_p,
\end{equation*}
\vspace{-0.35cm} 
	 \begin{equation*}
	 [t^{\bf{m}} d_i, t^{\bf{k}} d_j]= k_i t^{\bf{m+k}} d_j - m_jt^{\bf{m+k}} d_i + \phi(t^{\bf{m}} d_i, t^{\bf{k}} d_j),
\end{equation*}
	 \begin{equation*}
	[t^{\bf{m}}d_i, x \otimes t^{\bf{k}}]= k_i x \otimes t^{\bf{m+k}}, \quad \forall x \in \gg, \,\,\,\, {\bf{m}}, {\bf{k}} \in \Z^{n},\,\, 1 \leq i,j \leq n.
\end{equation*}
Let $B$ be a finitely generated commutative associative algebra with unit. In this paper our concerned Lie algebra is 
\begin{center}
	$\uptau(B):= \uptau \otimes B= (\mathcal{L}(\gg) \oplus  \mathcal{K} \oplus \mathcal{W}_n) \otimes B$,
\end{center}
called the map full toroidal Lie algebra with the natural bracket 
\begin{center}
	$[X \otimes a, Y \otimes b]=[X,Y] \otimes ab$, \quad $X,Y \in \uptau, \,\, a,b \in B$.
\end{center} 
For notational convenience, we denote $X \otimes a$ by $X(a)$. Unlike the full toroidal Lie algebra which has the finite dimensional center, $\uptau(B)$ has infinite dimensional center, given by the set Span $\{K_i(a): a \in B, \,1 \leq i \leq n\}$.\par
We now introduce an alternative set of notations for the elements of $\uptau$. For ${\bf{p}} \in \C^n, {\bf{m}} \in \Z^n$, set 
\begin{equation*}
	D({\bf{p}}, {\bf{m}})= \sum_{i=1}^{n} p_i t^{\bf{m}} d_i, \quad \quad K({\bf{p}}, {\bf{m}})= \sum_{i=1}^{n} p_i t^{\bf{m}} K_i.
\end{equation*}
 With these notations bracket formulas can be expressed as
\begin{align*}
	[D({\bf{p}}, {\bf{m}}), D({\bf{q}}, {\bf{k}})] &= D({\bf{w}}, {\bf{m+k}}) +\phi(D({\bf{p}}, {\bf{m}}), D({\bf{q}}, {\bf{k}})), \quad \text{where}\,\, {\bf{w}}=({\bf{p}}|{\bf{k}}){\bf{q}}-({\bf{q}}|{\bf{m}}){\bf{p}},\\
	\medskip
	[D({\bf{p}}, {\bf{m}}), K({\bf{q}}, {\bf{k}})] &= K({\bf{z}}, {\bf{m+k}}), \quad \text{where}\,\, {\bf{z}}=({\bf{p}}|{\bf{k}}){\bf{q}}+({\bf{q}}|{\bf{p}}){\bf{m}},\\
\vspace{0.3cm}
	[D({\bf{p}}, {\bf{m}}), x \otimes t^{\bf{k}}]&= ({\bf{p}}|{\bf{k}})\, x \otimes t^{\bf{m+k}}, \hspace{1cm} [K({\bf{p}}, {\bf{m}}), K({\bf{q}}, {\bf{k}})]=0
\end{align*}
Also, cocycles are given by
\begin{align*}
	\phi_1(D({\bf{p}}, {\bf{m}}), D({\bf{q}}, {\bf{k}}))&=-({\bf{p}}|{\bf{k}})({\bf{m}}|{\bf{q}})\, K({\bf{m}}, {\bf{m+k}}),	\medskip\\
	\phi_2(D({\bf{p}}, {\bf{m}}), D({\bf{q}}, {\bf{k}}))&=({\bf{p}}|{\bf{m}})({\bf{q}}|{\bf{k}})\, K({\bf{m}}, {\bf{m+k}}), \quad \forall {\bf{p}},{\bf{q}} \in \C^n,\, {\bf{m}},{\bf{k}} \in \Z^n.
\end{align*}
\par
Let $D= \text{Span}\, \{d_1, d_2, \ldots, d_n\}$, and define $\delta_i \in D^{*}$ such that $\delta_i(d_j)= \delta_{i,j}$ for $1 \leq i,j \leq n$. For ${\bf{m}}= (m_1,m_2, \ldots, m_n)$, set $\delta_{{\bf{m}}}= \sum_{i=1}^n m_{i}\, \delta_{i}$. For notational simplicity, we denote $\delta_{\mathbf{m}}$ by $\mathbf{m}$. Notice that $\uptau$ admits a natural $\Z^n$-grading:
\begin{equation}\label{hj}
\uptau= \displaystyle{\bigoplus_{{\bf{m}} \in \Z^n} \uptau_{{\bf{m}}}}, \quad \text{where} \quad \uptau_{{\bf{m}}}= \bigg( \gg \otimes t^{\bf{m}} \oplus \sum\limits_{i=1}^n \C t^{\textbf{m}} K_i  \oplus \hspace{0.1cm} \bigoplus\limits_{i=1}^n \C t^{\bf{m}}d_i \bigg).
\end{equation}

Similarly, $\uptau(B)$ also admits a natural $\Z^n$-grading:
\begin{equation}\label{ed}
	\uptau(B)= \displaystyle{\bigoplus_{{\bf{m}} \in \Z^n} \uptau_{{\bf{m}}}(B)}, \quad 
	\text{where } \uptau_{\mathbf{m}}(B) = \uptau_{\mathbf{m}} \otimes B.
\end{equation} 
Now we recall some basic definitions.
	\subsection{Trivial module} Let $L$ be a Lie algebra and let $V$ be an $L$-module. Then $V$ is called a trivial module if $x.v=0$ for all $x \in L$ and $v \in V$. A vector $v \in V$ is called a trivial vector if $x.v=0$ for all $x \in L$.
\subsection{Harish-Chandra Module} A $\uptau(B)$-module $V$ is said to be Harish-Chandra module if it satisfies the following:
\begin{enumerate}
	\item $V$ is a weight module with respect to $D= \text{Span}\, \{d_1, d_2, \ldots, d_n\}$, i.e.,
	\medskip
	\begin{center}
		$V= \oplus_{\lambda \in D^{*}} V_{\lambda}$, where $V_{\lambda}= \{v \in V: d_i.v=\lambda(d_i) v,\,\, 1 \leq i \leq n\}$.
	\end{center}
	\medskip
\item dim $V_{\lambda} < \infty$,\, for all ${\lambda} \in D^{*}$.
	\end{enumerate}
	The set $P(V)= \{\lambda \in D^{*}: V_{\lambda} \neq 0\}$ is called the set of weights of $V$ and any non-zero element $v \in V_{\lambda}$ is called a weight vector with weight $\lambda$. 
	\par
	\subsection{Cuspidal Module} A Harish-Chandra module $V$ is said to be cuspidal if the dimensions of its weight spaces are uniformly bounded, i.e., there exists $M \in \N$ such that dim $V_{\lambda} < M$, $\forall {\lambda} \in D^{*}$. These modules are also known as uniformly bounded modules.

	\subsection{Highest weight module} A $\uptau(B)$-module $V$ is called a highest weight module if there exists a nonzero vector $v \in V$ such that 
	\begin{center}
		$\uptau(B)^+v=0$ \,\,\,\, and \,\,\,\, $\mathcal{U}(\uptau(B))v=V$.
	\end{center}
	In this case, $v$ is called  a highest weight vector.
	\subsection{Generalized highest weight module} A $\uptau(B)$-module $V$ is called a generalized highest weight (GHW) module if there exists a nonzero vector $v \in V$ and $k \in \N$ such that 
	\begin{center}
		$\uptau_{{\bf{m}}}(B)v=0$ \,\,\,\, for all \, ${\bf{m}} \geq (k,k,\ldots,k)$.
	\end{center} 
	Such a vector $v$ is called a GHW vector.
	\subsection{Single point evaluation module}
	Let $V$ be a $\uptau(B)$-module with the corresponding representation $\rho: \uptau(B) \rightarrow \text{End}\,(V)$. Then $V$ is said to be a single point evaluation module for $\uptau(B)$ if there exists a maximal ideal $M$ of $B$ such that the map $\rho$ factors through $\uptau(M)$. In that case, $V$ can be viewed as a module over $\uptau$ as
	\begin{center}
		$\uptau(B)/ \uptau(M) \cong \uptau(B/M) \cong \uptau(\C) \cong \uptau$.
	\end{center}  Equivalently, if $V$ is a single point evaluation module for $\uptau(B)$, then there exists an algebra homomorphism $\phi: B \rightarrow \C$ such that
	\begin{center}
		$X(b).v=\phi(b)\,X.v$ \quad \text{for all} $X \in \uptau,\, b \in B$, and $v \in V$.
	\end{center}
	It is easy to check that $V$ is an irreducible single point evaluation $\uptau(B)$-module if and only if $V$ is an irreducible module for $\uptau$.

	\subsection{Change of Coordinates} Let $C = (c_{ij})$ be an $n \times n$ matrix in $GL_n{(\Z)}$ with determinant $\pm 1$. Then one has the following automorphism of $\uptau$, which we again denote by C:
	\begin{center}
		$C(x \otimes t^{\textbf{m}}) = x \otimes t^{\textbf{m}C^{T}}, \hspace{0.5cm}
		C(t^{\textbf{m}}K_i) = \displaystyle{\sum_{r = 1}^{n}{c_{r i}t^{\textbf{m}C^{T}}K_r}}, \hspace{0.5cm}
		C(t^{\textbf{m}}d_i) = \displaystyle{\sum_{r = 1}^n b_{ir}t^{\textbf{m}C^{T}}d_r},$
	\end{center}
	where $C^{-1} = (b_{ij})$ and $C^{T}$ denotes the transpose of $C$. This automorphism 
	$C$ yields a new full toroidal Lie algebra from the original one.  \par
	Alternatively, we can define another Lie algebra automorphism of $\uptau$ in a different way. Let $\{{\bm{\al}}_1,\bm{\al}_2,\ldots,\bm{\al}_n\}$ be a $\Z$-basis of $\mathbb{Z}^n$ and $\{\bm{\beta}_1,\bm{\beta}_2,\ldots, \bm{\beta}_n\}$ be the corresponding dual basis in $\mathbb{R}^n$. Define $\Phi: \uptau \rightarrow \uptau$ by
	\vspace{0.15cm}
	\begin{center}
		$x \otimes t^{\bf{m}} \mapsto x \otimes t^{\sum_{i=1}^{n}m_i\bm{\al}_i},\hspace{0.4cm} t^{\bf{m}}K_i \mapsto K(\bm{\al}_i,\, \sum_{i=1}^{n}m_i\bm{\al}_i) , \hspace{0.4cm} t^{\bf{m}}d_i \mapsto D(\bm{\beta}_i,\, \sum_{i=1}^{n}m_i\bm{\al}_i) $
	\end{center}
		\vspace{0.15cm}
The automorphisms $\Phi$ and $C$ of $\uptau$ naturally extend to automorphisms $\overline{\Phi}$ and  $\overline{C}$ of $\uptau(B)$, respectively. More explicitly, $\overline{\Phi}:\uptau(B) \rightarrow \uptau(B)$ is defined by \begin{center}
$\overline{\Phi}\,(X(b))= \Phi(X)(b),\quad \forall \,X \in \uptau,\, b \in B$.
\end{center} 
We call this process as a {\it change of co-ordinates}.
	\subsection{Representation of toroidal Lie algebras}
	Let $P_{\gg}^+$ denote the set of dominant integral weights of a Lie algebra $\gg$. For $\lambda_1 \in P_{\gg}^+$, let $V(\lambda_1)$ be the irreducible highest weight $\gg$-module with the highest weight $\lambda_1$. Similarly, let $\lambda_2 \in P_{\mathfrak{sl}_{n}}^+$, and $V(\lambda_2)$ be the highest weight module for ${\mathfrak{sl}_{n}}$. We extend $V(\lambda_2)$ to a ${\mathfrak{gl}_{n}}$-module by assuming the identity matrix acts by a scalar $c$; denote this module by $V(c, \lambda_2)$. Let ${\boldsymbol{\alpha}}=(\alpha_1,\ldots,\alpha_n) \in \C^n$. Define an action of $\uptau$ on the tensor product $V(\lambda_1) \otimes V(c, \lambda_2) \otimes \C[t_1^{\pm1}, t_2^{\pm1}, \ldots, t_n^{\pm1}]$ as follows
	\begin{center}
		$\big(x \otimes t^{\bf{m}}\big)(v_1 \otimes v_2 \otimes t^{\bf{r}}) = (x v_1) \otimes v_2 \otimes t^{\bf{m+r}}$, \hspace{1cm}$\big(t^{\bf{m}} K_i\big) (v_1 \otimes v_2 \otimes t^{\bf{r}}) =0$,\\
		\vspace{0.3cm}
		$\big(t^{\bf{m}} d_i\big) (v_1 \otimes v_2 \otimes t^{\bf{r}}) =(\alpha_i + r_i)(v_1 \otimes v_2 \otimes t^{\bf{m+r}}) + \sum_{j=1}^{n} m_j \,v_1 \otimes (E_{j,i} v_2) \otimes t^{\bf{m+r}},$\\
		\vspace{0.25cm}
		$\forall x \in \mathfrak{g},\, v_1 \in V(\lambda_1), v_2 \in V(c, \lambda_2), \,{\bf{m}},{\bf{r}} \in \mathbb{Z}^{n}, \,\,\, 1 \leq i \leq n.$
	\end{center}
	
	Here $ E_{j,i}$ denotes the $n \times n$ matrix with 1 in the $(j,i)$-th position and 0 elsewhere. We denote the above $\uptau$-module by $\mathfrak{L}(c, \lambda_1, \lambda_2, \boldsymbol{\alpha})$. 
	\par
	Proposition 3.17 of \cite{sou} shows that $\mathfrak{L}(c, \lambda_1, \lambda_2, \boldsymbol{\alpha})$ is an irreducible module for $\uptau$ if either $(\lambda_2, c, \boldsymbol{\alpha}) \notin \{0\} \times \{0,n\} \times \Z^n$ or $\lambda_1 \neq 0$ or $(\lambda_2,c) \neq (\omega_k, k)$ for any fundamental weight $\omega_k$ of $\mathfrak{sl}_n$ with $1 \leq k \leq n-1$.
	\par
	In particular, when $\lambda_1=0$, the module $\mathfrak{L}(c, 0, \lambda_2, \boldsymbol{\alpha})$ is an irreducible $\uptau$-module if and only if it is irreducible as a module over $\mathcal{W}_n$. In the literature, these modules are known as Larson modules or Shen modules \cite{ta,sh}. Let $\{\omega_1, \omega_2,\ldots, \omega_{n-1}\}$ be the set of fundamental weights of $\mathfrak{sl}_n$. Also, we assume that $\omega_0=\omega_n=0$ for notational convenience. Let \( W \cong \mathbb{C}^{n} \) denote the standard representation of \( \mathfrak{gl}_{n} \). Then its $k$-fold exterior power $\wedge^k W$ is known to be isomorphic to the irreducible $\mathfrak{gl}_{n}$-module $V(k,\omega_k)$ for each $0 \leq k \leq n$. The associated $\mathcal{W}_n$-modules $\mathfrak{L}(k, 0, \omega_k, \boldsymbol{\alpha})$ consist of differential $k$-forms and they form the de Rham complex
	\begin{center}
		$\mathfrak{L}(0, 0, \omega_0, \boldsymbol{\alpha}) \xrightarrow{d_0} \mathfrak{L}(1, 0, \omega_1, \boldsymbol{\alpha}) \xrightarrow{d_1} \cdots \xrightarrow{d_{n-1}} \mathfrak{L}(n, 0, \omega_{n}, \boldsymbol{\alpha})$
	\end{center}
	
	The differential map $d_k:\mathfrak{L}(k, 0, \omega_k, \boldsymbol{\alpha}) \rightarrow \mathfrak{L}(k+1, 0, \omega_{k+1}, \boldsymbol{\alpha})$ is a homomorphism of $\mathcal{W}_n$-modules and consequently the kernel and image of the map naturally give rise to $\mathcal{W}_n$-modules of $\mathfrak{L}(k, 0, \omega_k, \boldsymbol{\alpha})$ and $\mathfrak{L}(k+1, 0, \omega_{k+1}, \boldsymbol{\alpha})$, respectively. By the Proposition 3.17 of \cite{sou}, we know that for any $0 \leq k \leq n$, the module $\mathfrak{L}(k, 0, \omega_k, \boldsymbol{\alpha})$ admits a unique irreducible quotient over $\mathcal{W}_n$. This irreducible quotient is isomorphic to the image of the map $d_k:\mathfrak{L}(k, 0, \omega_k, \boldsymbol{\alpha}) \rightarrow \mathfrak{L}(k+1, 0, \omega_{k+1}, \boldsymbol{\alpha})$ and is denoted by $d_k\mathfrak{L}(k, 0, \omega_k, \boldsymbol{\alpha})$. Furthermore, if $\boldsymbol{\alpha} \in \Z^n$, then the module $\mathfrak{L}(n, 0, \omega_{n}, \boldsymbol{\alpha})$ admits a one dimensional trivial module as its irreducible quotient. We now state the following classification theorem as given in Theorem 3.24 of \cite{sou}.
	
	\begin{thm}[\cite{sou}, Theorem 3.24]\label{fg}
		Let $V$ be a non trivial irreducible cuspidal $\uptau$-module. Then $V$ is isomorphic to one of the following:
		\begin{enumerate}
			\item $\mathfrak{L}(c, \lambda_1, \lambda_2, \boldsymbol{\alpha})$, for some quadruplet $(c, \lambda_1, \lambda_2, \boldsymbol{\alpha}) \in \C \times (P^+_{\gg})^{\times} \times (P^+_{\mathfrak{sl}_n}) \times \C^n$.
			\item $\mathfrak{L}(c, 0, \lambda_2, \boldsymbol{\alpha})$ for some $(c, \lambda_2, \boldsymbol{\alpha}) \in \C \times (P^+_{\mathfrak{sl}_n}) \times \C^n$ with $(c, \lambda_2) \neq (k, \omega_k)$ for $1 \leq k \leq n-1$.
			\item $\mathfrak{L}(n, 0, 0, \boldsymbol{\alpha})$ with $\boldsymbol{\alpha} \notin \Z^n$.
			\item A submodule $d_k\mathfrak{L}(k, 0, \omega_k, \boldsymbol{\alpha})$ of $\mathfrak{L}(k+1, 0, \omega_{k+1}, \boldsymbol{\alpha})$ for some $\boldsymbol{\alpha} \in \C^n$, where $0 \leq k \leq n-1$ and $\omega_0=\omega_n=0$.
		\end{enumerate}
	\end{thm}
	
	\section{Classification of Irreducible cuspidal modules over $\uptau(B)$}
	In this section, we classify all the irreducible cuspidal modules over $\uptau(B)$. Note that if $V$ is an irreducible $\uptau(B)$-module, then there exists $\lambda \in D^{*}$ such that $P(V) \subset \{\lambda+ {\bf{m}}: {\bf{m}} \in \Z^{n}\}$. Throughout this section, $V$ always stands for an irreducible cuspidal $\uptau(B)$-module.
	
	\begin{prop}\label{po}
		$\mathcal{K}(B)$ acts trivially on $V$.
	\end{prop}
	\begin{proof}
		First, we prove that zero degree central operators $K_i(B)$ act trivially on $V$ for each $1 \leq i \leq n$. Fix $1 \leq i \leq n$, $a \in B$, and choose  $h_1, h_2 \in \mathfrak{h}$ such that $\langle h_1, h_2 \rangle \neq 0$. Consider the Heisenberg algebra $\mathcal{H}(a)$ spanned by the elements $\{h_1 \otimes t_i^{n}, h_2 \otimes t_i^{m}(a), K_i(a): n>0, m<0, n,m \in \Z\}$, with the Lie bracket 
		\begin{center}
			$[h_1 \otimes t_i^{n}, h_2 \otimes t_i^{m}(a)]= \langle h_1, h_2 \rangle \delta_{n,-m} K_i(a)$, 
		\end{center}
		 and $K_i(a)$ acting centrally. Define $V'= \oplus_{m \in \Z} V_{\lambda + m{\bf{e}}_i}$. Then $V'$ is a $\Z$-graded $\mathcal{H}(a)$-module with uniformly bounded weight spaces. If $K_i(a)$ acts by a nonzero scalar on $V'$, then by Proposition 4.3(i) and 4.5 of \cite{F}, $V'$ must contain a highest weight or lowest weight module, contradicting the fact that $V'$ is uniformly bounded. Hence $K_i(a)$ acts trivially on $V'$ and thus on $V$, since $V$ is irreducible. \par
		\vspace{0.2cm}
		
		Next, Lemmas 3.6, 3.7, 3.8 of \cite{sou} show that $t^{\bf{p}}K_i$ acts locally nilpotently on $V$ for all $1 \leq i \leq n$ and $\bf{p} \in \Z^n$ with $p_j \neq 0$ for some $j \neq i$. Fix $1 \leq i \leq n$ and choose $\bf{p} \in \Z^n$ such that $p_j \neq 0$ for $1 \leq j \neq i \leq n$. As $V$ is a uniformly bounded module, there exists an integer $N \in \N$ such that 
		\begin{center}
			$((t^{\bf{-p}}K_i)(t^{\bf{p}}K_i)^N)V=\{0\}$.
		\end{center} 
		Then for any $a \in B, {\bf{m}} \in \Z^n$, and $v \in V$,
		
		\begin{equation*}
		\begin{aligned}
			0&=[t^{\bf{m+p}} \,d_j(a), ((t^{\bf{-p}}K_i)(t^{\bf{p}}K_i)^N)]v\\
			&=Np_j\bigg((t^{\bf{-p}}K_i)^N t^{\bf{m+2p}} K_i(a)(t^{\bf{p}}K_i)^{N-1} -t^{\bf{m}} K_i(a) (t^{\bf{-p}}K_i)^{N-1} (t^{\bf{p}}K_i)^{N}\bigg)v
		\end{aligned}
		\end{equation*}
		
		Applying $t^{\bf{p}}K_i$ to both sides, we get
		\begin{center}
			$(t^{\bf{m}} K_i(a)) \,(t^{\bf{-p}}K_i)^{N-1} (t^{\bf{p}}K_i)^{N+1}V=\{0\}$,\quad $\forall a \in B, {\bf{m}} \in \Z^n$.
		\end{center} 
 
 Therefore, for $a_1, a_2 \in B,\, {\bf{m}}^{(1)}, {\bf{m}}^{(2)} \in \Z^n$, and $v \in V$, we have
	
		\begin{equation*}
	\begin{aligned}
	0&=[t^{\bf{m^{(2)} +p}} \,d_j(a_2), (t^{\bf{m^{(1)}}} K_i(a_1)) \,(t^{\bf{-p}}K_i)^{N-1} (t^{\bf{p}}K_i)^{N+1}]v\\
	&=(N+1)\,p_j\, t^{\bf{m^{(1)}}} K_i(a_1)\, t^{\bf{m^{(2)}}+2\bf{p}} K_i(a_2)\, (t^{\bf{-p}}K_i)^{N-1} (t^{\bf{p}}K_i)^N v\\
	&\,\,\,\,\,\, -(N-1)\,p_j\, t^{\bf{m^{(1)}}} K_i(a_1)\, t^{\bf{m^{(2)}}} K_i(a_2)\, (t^{\bf{-p}}K_i)^{N-2} (t^{\bf{p}}K_i)^{N+1} v.
	\end{aligned}
	\end{equation*}
	
	Applying $t^{\bf{p}}K_i$ again on both sides, we obtain 
	\begin{center}
		$(t^{\bf{m}^{(1)}} K_i(a_1))(t^{\bf{m^{(2)}}} K_i(a_2)) \,(t^{\bf{-p}}K_i)^{N-2} (t^{\bf{p}}K_i)^{N+2}V=\{0\}$, \hspace{0.25cm} $\forall a_1, a_2 \in B, {\bf{m}^{(1)}}, {\bf{m}^{(2)}} \in \Z^n$.
	\end{center} 
	 By continuing this process $(N-2)$ times, we deduce 
	 
	 \begin{equation}\label{gh}
	 	(t^{\bf{m}^{(1)}} K_i(a_1)) \cdots (t^{\bf{m^{(N)}}} K_i(a_N)) (t^{\bf{p}}K_i)^{2N}V=\{0\}, \quad \forall a_k \in B, {\bf{m}^{(k)}} \in \Z^n.
	 \end{equation} 
	  Now by applying successively $t^{\bf{m}^{(N+1)}-\bf{p}}\, d_j(a_{N+1}), t^{\bf{m}^{(N+2)}-\bf{p}}\, d_j(a_{N+2}),\ldots, t^{\bf{m}^{(3N)}-\bf{p}}\, d_j(a_{3N})$ to both sides of \eqref{gh} and following the same method as above, we obtain
	
	\begin{equation}\label{gk}
	(t^{\bf{m}^{(1)}} K_i(a_1))(t^{\bf{m^{(2)}}} K_i(a_2)) \cdots (t^{\bf{m^{(3N)}}} K_i(a_{3N}))V=\{0\}, \quad \forall a_k \in B,\, {\bf{m}^{(k)}} \in \Z^n.
	\end{equation}
Since $i$ is arbitrary, equation \eqref{gk} holds for every $1 \leq i \leq n$. Using arguments similar to those in \cite[Lemma 4.4]{sp}, we can find a nonzero vector $v \in V$ such that $(t^{\bf{m}} K_i(a))v=0$ for all $a \in B, {\bf{m}} \in \Z^n$, and $1 \leq i \leq n$. Then $W=\{v \in V: \mathcal{K}(B)v=0 \}$ becomes a nonzero $\uptau(B)$-submodule of $V$. Hence by the irreducibility of $V$, the proof follows.
	\end{proof}
	Our next aim is to prove the existence of a maximal ideal $\mathcal{M}$ of $B$ such that $\uptau(\mathcal{M})$ acts trivially on $V$. The following proposition is the first step towards that direction.
	\begin{prop}\label{va}
		There exists a cofinite ideal $J$ of $B$ such that $(\uptau \otimes J) V=0$.
	\end{prop}

\begin{proof}
	If $V$ is trivial, there is nothing to prove. So assume $V \neq 0$. Consider the weight space decomposition of $V= \bigoplus_{{\bf{m}} \in \Z^{n}} V_{\lambda + \bf{m}}$, for some $\lambda \in D^*$. Choose ${\bf{m}} \in \Z^n$ such that $V_{\lambda + \bf{m}} \neq 0$. For $j \in \Z \setminus \{0\}$, define
	\begin{center}
		$I_{j}= \{a \in B: (t_1^{j}d_i(a)) V_{\lambda + \bf{m}}=0, \,\,\forall \,1 \leq i \leq n\}$.
	\end{center}
\noindent\textbf{Claim (a).} $I_{j}$ is a cofinite ideal of $B$.\\
To see that $I_{j}$ is an ideal, take $a \in I_{j}, b \in B$, and $v \in V_{\lambda + \bf{m}}$. Then
\begin{center}
	$t_1^{j} d_i(ab)v=\frac{1}{j} [d_1(b), t_1^{j} d_i(a)]v=0, \,\,\,\,\, \forall\, 1 \leq i \leq n$ .
\end{center}
In the above expression we have used the fact that $d_1(b)$ preserves weight spaces. Note that $I_j$ is the kernel of the linear map
\begin{center}
	$\Phi: B \rightarrow (Hom_{\C}\, (V_{\lambda + \bf{m}}, V_{\lambda + {\bf{m}}+ j{\bf{e}}_1}))^n$, \quad $\Phi(a)= (f_1(a), f_2(a), \ldots ,f_n(a))$,
\end{center} 
where $f_i(a)(v)= (t_1^{j}d_i(a))v$, \,$\forall \, 1 \leq i \leq n$.
Since $V$ is uniformly bounded, there exists $N \in \N$ such that $\text{dim}\,V_{\lambda + \bf{m}} \leq N$ for all ${\bf{m}} \in \Z^{n}$. Thus we have dim $B/I_j \leq N^{2n}$ and that proves our claim.\par
\noindent\textbf{Claim (b).} $I_1^{j}I_2 \subset I_{j+2}$ for all $j \geq 1$.\\ We prove this by induction on $j$. For $j=1$, we need to show that $a \in I_1$ and $b \in I_2$ implies $ab \in I_3$. Indeed, this holds because
\begin{center}
	$0=[t_1d_1(a), t_1^{2}d_i(b)]V_{\lambda + \bf{m}}= (2t_1^{3}d_i(ab) - \delta_{1,i}t_1^{3}d_1(ab))V_{\lambda + \bf{m}}, \quad \forall \,1 \leq i \leq n$
\end{center}
Now assume that the claim holds for some fixed $j \geq 1$. Let $a \in I_1, b \in I_1^{j}I_2$, then $t_1^{j+2}d_i(b)v=0$ for all $i$, as $I_1^{j}I_2 \subset I_{j+2}$. Then by the induction the claim follows since
\begin{center}
	$0=[t_1d_1(a),\, t_1^{j+2}d_i(b)]V_{\lambda + \bf{m}}= (j+2)\,t_1^{j+3}d_i(ab)V_{\lambda + \bf{m}} - \delta_{1,i}\, t_1^{j+3}d_1(ab)V_{\lambda + \bf{m}}$.
\end{center}
Now by using an argument exactly similar to that of Proposition 4.1 of \cite{SAV}, we have 
$I_1^{N^{2n}}I_2 \subset I_{j}$ for $j \geq 1$. Proceeding similarly as above, we get $I_{-1}^{N^{2n}}I_{-2} \subset I_{-j}$ for $j \geq 1$. Take $J_1= I_{-1}^{N^{2n}}I_{-2}I_1^{N^{2n}}I_2$, then $J_1 \subset I_j$ for all $j \neq 0$. Since $J_1$ is the product of finitely many cofinite ideals, it is cofinite. Now any element of $J_1$ can be written as a sum of the elements of the form $a_{-1}a_{1}$ with $a_{-1} \in I_{-1}, a_1 \in I_1$ since $J_1 \subset I_{-1}I_1$. For $1 \leq i \leq n$, we get 
\begin{center}
	$0=[t_1^{-1}d_1(a_{-1}), t_1d_i(a_1)] V_{\lambda + \bf{m}}= (d_i(a_{-1}a_{1}) + \delta_{1,i}d_1(a_{-1}a_{1})) V_{\lambda + \bf{m}}$.
\end{center}
Thus $d_i(J_1) V_{\lambda + \bf{m}}=0$ for all $1 \leq i \leq n$.
 Also, since $J_1 \subset I_j$ for all $j \neq 0$, we have $t_1^{j}d_i(J_1) V_{\lambda + \bf{m}}=0$ for all $j \neq 0$ and $1 \leq i \leq n$. Similarly we can construct cofinite ideals $J_2, J_3, \ldots , J_n$ of $B$ such that
  \begin{center}
 	$t_k^{j}\, d_i(J_k)V_{\lambda + \bf{m}}=0$ \quad for all $j \in \Z,\,\, 2 \leq k \leq n,\,\, 1 \leq i \leq n$.
 \end{center}
 Since as a Lie algebra $\mathcal{W}_n$ is generated by the set $\{t_k^{j}\, d_i: j \in \Z,\, 1 \leq i,k \leq n\}$, we conclude that
 \begin{center}
 	$\mathcal{W}_n(J')\,V_{\lambda + \bf{m}}=0$,\quad where $J'=J_1J_2\cdots J_n$.
 \end{center}  Next, consider $J''=\{a \in B: (x \otimes t_1)(a)V_{\lambda + \bf{m}}=0, \,\forall x \in \gg \}$, then it is straightforward to verify that $J''$ is a cofinite ideal of $B$. Let $a \in J', b \in J''$. Then \begin{center}
	$(x \otimes t^{\bf{m}})(ab) V_{\lambda + \bf{m}} =[t^{{\bf{m}}- {\bf{e}}_1}\,d_1(a), x \otimes t_1(b)] V_{\lambda + \bf{m}}=0$.
\end{center}
Therefore we get $\uptau(J) V_{\lambda + \bf{m}}=0$ by taking $J=J'J''$.\par
Since $V_{\lambda + \bf{m}} \neq 0$, and $V$ is irreducible, we have $V= \mathcal{U}(\uptau(B)) V_{\lambda + \bf{m}}$. Let $\mathcal{U}_0(\uptau(B)) \subseteq \mathcal{U}_1(\uptau(B)) \subseteq \mathcal{U}_2(\uptau(B)) \subseteq \cdots$ denote the natural filtration of $\uptau(B)$, induced by the grading on the tensor algebra of $\uptau(B)$. To prove $\uptau(J)V=0$, it is enough to prove that $\uptau(J)\, \mathcal{U}_k(\uptau(B)) V_{\lambda + \bf{m}}=0,\,\, \forall k \in \N$. An easy induction on $k$ completes the proof.
	\end{proof}
	
Now, since $J$ is a cofinite ideal, the quotient ring $B/J$ is an Artinian ring. Hence $B/J$ can be written as a finite product of local Artinian rings. It follows that the ideal $J$ is contained in a finite product of powers of maximal ideals. Therefore, we may assume that $J= \mathcal{M}_1^{k_1} \mathcal{M}_2^{k_2} \cdots \mathcal{M}_r^{k_r}$, where each $\mathcal{M}_i$ is a maximal ideal of $B$ and $k_i \in \N$. We have the following proposition:
\begin{prop}\label{vb}
	Let $\M_1, \M_2$ be two coprime ideals of $B$ such that $\uptau(\M_1 \M_2) V=0$, then either $\uptau(\M_1)V=0$ or $\uptau(\M_2)V=0$.
\end{prop}
\begin{proof}
	Since $\uptau(\M_1 \M_2) V=0$, we may view $V$ as a module over $\uptau(B/\M_1 \M_2)$. As $\M_1$ and $\M_2$ are coprime, the Chinese Remainder Theorem gives
	\begin{center}
		$\uptau(B/\M_1 \M_2) \cong \uptau(B/\M_1) \oplus \uptau(B/\M_2)$.
	\end{center} Thus $V$ becomes a module over $\uptau(B/\M_1) \oplus \uptau(B/\M_2)$ via this isomorphism. Let $V= \oplus_{{\bf{m}} \in \Z^n} V_{\lambda+ \bf{m}}$ be the weight space decomposition of $V$. For $1 \leq i \leq n$, define 
	\begin{center}
		$d_i^1= d_i \otimes (1+\M_1),\quad d_i^2= d_i \otimes (1+\M_2),$
	\end{center} 
	so that $d_i= d_i^1 + d_i^2$. By an argument similar to that in claim 1 of Theorem 4.3 of \cite{ZRX}, both $d_i^1 , d_i^2$ act diagonalizably on $V$. As $d_i^1 , d_i^2$ commute with each other, we have the following decomposition of $V= \oplus_{\bf{p}, \bf{q} \in \Z^n} V_{\bf{p}, \bf{q}}$, where
	\begin{center}
		$V_{\bf{p}, \bf{q}}= \{v \in V: d_i^1. v= (\lambda_1(d_i)+ p_i)v, \,\, d_i^2. v= (\lambda_2(d_i)+ q_i)v, \,\, 1 \leq i \leq n\}$ .
	\end{center}
	with $\lambda_1 + \lambda_2= \lambda$. Observe that $V_{\lambda + \bf{m}}= \oplus_{\bf{p} \in \Z^n} V_{\bf{m}+ \bf{p}, - \bf{p}}$. Define the subspaces $V_{(\bf{p})}$ and $V^{(\bf{q})}$ by
	\begin{equation*}
		V_{(\bf{p})}= \bigoplus_{\bf{q} \in \Z^n} V_{\bf{p}, \bf{q}} \hspace{0.5cm} and \hspace{0.5cm} V^{(\bf{q})}= \bigoplus_{\bf{p} \in \Z^n} V_{\bf{p}, \bf{q}}.
	\end{equation*}
	Then it is easy to see that $V^{(\bf{q})}$ and $V_{(\bf{p})}$ are uniformly bounded modules for $\uptau(B/\M_1)$ and $\uptau(B/\M_2)$ respectively.  Note that our claim is equivalent to proving that either $\uptau(B/\M_1)$ or $\uptau(B/\M_2)$ acts trivially on $V$. Assume, for contradiction, that this is not the case. \par 
	Let $\uptau(B/\M_1) V^{(\bf{q})}=0$ for some ${\bf{q}} \in \Z^{n}$, then by using the irreducibility of $V$, we have $V= \U(\uptau(B/\M_2)) V^{(\bf{q})}$. Since $\uptau(B/\M_1)$ commutes with $\uptau(B/\M_2)$, it follows that $\uptau(B/\M_1)$ acts trivially on $V$. This contradicts our assumption, so we must have $\uptau(B/\M_1) V^{(\bf{q})} \neq 0$ for all ${\bf{q}} \in \Z^n$. Now choose ${\bf{q}} \in \Z^n$ such that $V^{(\bf{q})} \neq 0$. Then $V^{(\bf{q})}$ is a nonzero uniformly bounded $\uptau(B/\M_1)$ module. Since $\uptau(B/\M_1)$ contains a copy of the Virasoro algebra, Lemma 2.11 and Corollary 2.13 of \cite{SAV}, imply that $V_{\bf{p}, \bf{q}} \neq 0$ for ${\boldsymbol{\lambda_1}} + \bf{p} \neq 0$, where $\boldsymbol{\lambda_1}=(\lambda_1(d_1),\ldots,\lambda_1(d_n))$. This implies that if ${\boldsymbol{\lambda_1}} + \bf{p} \neq 0$ then $V_{({\bf{p}})} \neq {\bf{0}}$. Again, using the same arguments, we get
	\begin{center}
		$V_{\bf{p}, \bf{q}} \neq 0$ \quad whenever ${\boldsymbol{\lambda_1}} + \bf{p} \neq 0$ and ${\boldsymbol{\lambda_2}} + \bf{q} \neq 0$.
	\end{center} 
	But this contradicts the fact that $V$ is a Harish-Chandra module since $V_{\lambda}= \oplus_{\bf{p} \in \Z^n} V_{\bf{p}, \bf{-p}}$ would then be infinite dimensional. This contradiction proves the proposition.
	\end{proof}	 

By the above proposition, we may assume that there exists a maximal ideal $\mathcal{M}$ and an integer $k \in \N$ such that $\uptau(\mathcal{M}^k) V=0$. The next proposition ensures the existence of the desired maximal ideal.
\begin{prop}\label{vc}
	Let $\mathcal{M}$ be an ideal of $B$ such that $\uptau(\mathcal{M}^{k}) V=0$ for some $k \in \N$, then we have $\uptau(\mathcal{M}) V=0$.
\end{prop}
\begin{proof}
	We prove this statement by using induction on $k$. For $k=1$, the result is obvious. Assume $k=2$, then the elements of 
	$\uptau(\M)$ act as commuting operators on $V$. Note that the operator $d_i(b)$ fixes the weight spaces of $V$. Then, using the same technique used in Proposition 4.5 of \cite{SAV}, we obtain
	\begin{center}
		$(d_i(b))^N\, V=0, \quad \forall \, 1 \leq i \leq n,\, b \in \M$,
	\end{center}
	where $N= \text{max} \, \{\,\text{dim}\, V_{\lambda+ \bf{m}}: {\bf{m}} \in \Z^n \}$.
	
	\par
	\vspace{0.2cm}
	\noindent\textbf{Claim 1.} For all $1 \leq i \leq n$ and $p \in \Z$, we have $(t_i^{p}d_i(\M))^N V=0$.
	
	For $p=0$, this is already shown. Let $p \neq 0$ and take $b \in \M, v \in V$. Then, 
	\begin{center}
		$0=t_i^{p}d_i(d_i(b))^N v=[t_i^{p}d_i, (d_i(b))^N ] v= -Np  t_i^{p}d_i(b)(d_i(b))^{N-1} v$.
	\end{center}
Thus, $t_i^{p}d_i(b)(d_i(b))^{N-1} V=0$, and by applying $t_i^{p}d_i$ on $t_i^{p}d_i(b)(d_i(b))^{N-1} V=0$ repeatedly $(N-1)$ more times, the claim follows.\par
	\vspace{0.2cm}
\noindent\textbf{Claim 2.}  For $1 \leq i \leq n$, $0 \leq r \leq N$ and $a \in \M$, we have
\begin{center}
	$(t^{{\bf{m}}^{(1)}}d_i(a))(t^{{\bf{m}}^{(2)}}d_i(a))\cdots (t^{{\bf{m}}^{(r)}}d_i(a))(t_i^{p}d_i(a))^{N-r} V=0$,
\end{center}  
holds for  all ${\bf{m}}^{(j)}=(m_1^{(j)},m_2^{(j)},\ldots, m_n^{(j)}) \in \Z^{n}$ and $p > \text{max}\, \{m_{i}^{(j)}: 1 \leq j \leq r\}$.
\par
For $r=0$, this is exactly Claim 1. Now assume that the result is true for some fixed $0 \leq r < N$. Then

\begin{equation*}
	\begin{aligned}
		0&=[t^{({\bf{m}}^{(r+1)}- p{\bf{e}}_i)}d_i, (t^{{\bf{m}}^{(1)}}d_i(a))(t^{{\bf{m}}^{(2)}}d_i(a))\cdots (t^{{\bf{m}}^{(r)}}d_i(a))(t_i^{p}d_i(a))^{N-r}]V\\
		&=\sum_{j=1}^{r} [t^{({\bf{m}}^{(r+1)}- p{\bf{e}}_i)}d_i, t^{{\bf{m}}^{(j)}}d_i(a)] t^{{\bf{m}}^{(1)}}\hspace{-0.05cm}d_i(a)\cdots t^{{\bf{m}}^{(j-1)}}\hspace{-0.05cm}d_i(a)t^{{\bf{m}}^{(j+1)}}\hspace{-0.05cm}d_i(a) \cdots t^{{\bf{m}}^{(r)}}\hspace{-0.05cm}d_i(a)(t_i^{p}d_i(a))^{N-r}V\\
		&\hspace{3cm}+ (t^{{\bf{m}}^{(1)}}d_i(a))(t^{{\bf{m}}^{(2)}}d_i(a))\cdots (t^{{\bf{m}}^{(r)}}d_i(a)) [t^{({\bf{m}}^{(r+1)}- p{\bf{e}}_i)}d_i, (t_i^{p}d_i(a))^{N-r}] V\\
	\end{aligned}
\end{equation*}
\begin{equation*}
\begin{aligned}
	&=\sum_{j=1}^{r} (m_i^{(j)}+p-m_i^{(r+1)}) \big(t^{({\bf{m}}^{(r+1)}+{\bf{m}}^{(j)}- p{\bf{e}}_i)}d_i(a)\big) t^{{\bf{m}}^{(1)}}\hspace{-0.05cm}d_i(a)\\
	 &\hspace{5cm}\cdots t^{{\bf{m}}^{(j-1)}}d_i(a)t^{{\bf{m}}^{(j+1)}}d_i(a) \cdots t^{{\bf{m}}^{(r)}}d_i(a)(t_i^{p}d_i(a))^{N-r}\big)V\\
	&\hspace{0.5cm}+(t^{{\bf{m}}^{(1)}}d_i(a))(t^{{\bf{m}}^{(2)}}d_i(a))\cdots (t^{{\bf{m}}^{(r)}}d_i(a)) (N-r)(2p-m_i^{(r+1)})t^{{\bf{m}}^{(r+1)}}d_i(a) (t_i^{p}d_i(a))^{N-r-1}V\\
	 &= (N-r)(2p-m_i^{(r+1)})(t^{{\bf{m}}^{(1)}}d_i(a))(t^{{\bf{m}}^{(2)}}d_i(a))\cdots (t^{{\bf{m}}^{(r)}}d_i(a)) t^{{\bf{m}}^{(r+1)}}d_i(a) (t_i^{p}d_i(a))^{N-r-1}V
	\end{aligned}
	\vspace{0.25cm}
	 \end{equation*}
	 Here in the third equality, each term of the summation is zero by induction hypothesis as $(m_i^{(r+1)}+ m_i^{(j)}-p) <p$ for all $j$. So the inductive step completes since $2p-m_i^{(r+1)} \neq 0$. Now by choosing some suitable $p$, we can deduce that 
 \begin{equation}\label{ad}
 	(t^{{\bf{m}}^{(1)}}d_i(a))(t^{{\bf{m}}^{(2)}}d_i(a))\cdots (t^{{\bf{m}}^{(N)}}d_i(a)) V=0, \,\,\,\,\forall\, 1 \leq i \leq n,\, {\bf{m}}^{(j)} \in \Z^n,\, a \in \M.
 \end{equation}
Let $R'=n(N-1)+1$, then by using \eqref{ad} we get
\begin{equation}\label{aa}
	(t^{{\bf{m}}^{(1)}}d_{i_1}(a))\cdots (t^{{\bf{m}}^{(M')}}d_{i_{R'}}(a)) V=0 \quad \text{for} \,\,\,a \in \M,\, 1 \leq i_j \leq n, \,{\bf{m}}^{(j)} \in \Z^{n}.
\end{equation}
Now using the same method as in claim 2, we deduce
\begin{equation}\label{ae}
	 (x_{1} \otimes t^{{\bf{m}}^{(1)}}(a))\cdots (x_{N} \otimes t^{{\bf{m}}^{(N)}}(a))V=0, \quad \forall x_{j} \in \gg,\, {\bf{m}}^{(j)} \in \Z^{n},\, a \in \M.
	\end{equation}
Let $R=R'+N$. Combining \eqref{aa} and \eqref{ae}, we obtain
\begin{equation}\label{ab}
	(x_{1} \otimes t^{{\bf{m}}^{(1)}}(a))\cdots (x_{r} \otimes t^{{\bf{m}}^{(r)}}(a))(t^{{\bf{m}}^{(r+1)}}d_{i_{r+1}}(a))\cdots (t^{{\bf{m}}^{(R)}}d_{i_R}(a)) V=0 ,
\end{equation} 
for all $a \in \M,\, 1 \leq i_j \leq n,\, {\bf{m}}^{(j)} \in \Z^{n},\, x_j \in \gg$. \par
\vspace{0.2cm}
\noindent\textbf{Claim 3.} There exists $P \in \N$ such that $(\uptau({\M}))^P\,\, V=0$.\vspace{0.1cm}\par
 Since the ideal $\mathcal{M}$ is cofinite, $\mathcal{M}^2$ is also cofinte. It follows that $\mathcal{M}/\mathcal{M}^2$ is finite dimensional as $\text{dim}\, (\mathcal{M}/\mathcal{M}^2) \leq \text{dim}\,(B/\mathcal{M}^2)$. Let $P= \text{dim}\, (\M/ \M^2)(R-1)+1$, then by expanding in a basis for $\M$ and using \eqref{ab}, we deduce that
\begin{equation*}
	(x_{1} \otimes t^{{\bf{m}}^{(1)}}(a_1))\cdots (x_{r} \otimes t^{{\bf{m}}^{(r)}}(a_r))(t^{{\bf{m}}^{(r+1)}}d_{i_{r+1}}(a_{r+1}))\cdots (t^{{\bf{m}}^{(P)}}d_{i_P}(a_P)) V=0 ,
\end{equation*}
for $a_j \in \M,\, 1 \leq i_j \leq n,\, {\bf{m}}^{(j)} \in \Z^{n},\, x_j \in \gg$. Hence the claim follows.\vspace{0.2cm}\par
Since $\uptau({\M})$ is an ideal of $\uptau(B)$, it is easy to see that 
\begin{equation}\label{ac}
	 \Big(\U(\uptau(B))\,(\uptau({\M}))\,\U(\uptau(B))\Big)^P=\U(\uptau(B))\,(\uptau({\M}))^P\,\U(\uptau(B)). 
	\end{equation}
Note that $(\U(\uptau(B))\,(\uptau({\M}))\,\U(\uptau(B)))\,V$ is a $\uptau(B)$-submodule of V and by \eqref{ac} it follows that $(\U(\uptau(B))\,(\uptau({\M}))\,\U(\uptau(B))) V \neq V$. Hence $(\U(\uptau(B))\,(\uptau({\M}))\,\U(\uptau(B)))\, V=0$ and that implies $\uptau({\M}) V=0$ as desired.\par
\vspace{0.2cm}
So the case $k=2$ is done, now assume the result is true for any $r<k$. Since $\uptau(\M^k)V=0$, we can consider $V$ as a module for $\uptau(B/ \M^k)$. Consider the ideal $\M^{k-1}/ \M^k$ of $B/ \M^k$, then we have $\uptau((\M^{k-1}/ \M^k)^2)V=0$ as $k \geq 2$. Then from the case $k=2$, we again have $\uptau(\M^{k-1}/ \M^k)V=0$. Hence $\uptau(\M^{k-1})V=0$, and by induction, it follows that $\uptau(\M)V=0$.
\end{proof}

Now, by the previous proposition, we may view $V$ as a module over $\uptau(A/\M) \cong \uptau$. Hence $V$ becomes a single point evaluation module. The following theorem is the main result of this section:
\begin{theorem}\label{ds}
	Let $V$ be an irreducible cuspidal module for $\uptau(B)$. Then $V$ is a single point evaluation module and is isomorphic to one of the following:
	\begin{enumerate}
		\item $\mathfrak{L}(c, \lambda_1, \lambda_2, \boldsymbol{\alpha}, \psi)$, for some quadruplet $(c, \lambda_1, \lambda_2, \boldsymbol{\alpha}) \in \C \times (P^+_{\gg})^{\times} \times (P^+_{\mathfrak{sl}_n}) \times \C^n$.
		\medskip
		\item $\mathfrak{L}(c, 0, \lambda_2, \boldsymbol{\alpha}, \psi)$ for some $(c, \lambda_2, \alpha) \in \C \times (P^+_{\mathfrak{sl}_n}) \times \C^n$ with $(c, \lambda_2) \neq (k, \omega_k)$ for $1 \leq k \leq n-1$.
		\medskip
		\item $\mathfrak{L}(n, 0, 0, \boldsymbol{\alpha}, \psi)$ with $\boldsymbol{\alpha} \notin \Z^n$.
		\medskip
		\item A submodule $d_k\mathfrak{L}(k, 0, \omega_k, \boldsymbol{\alpha}, \psi)$ of $\mathfrak{L}(k+1, 0, \omega_{k+1}, \boldsymbol{\alpha}, \psi)$ for some $\boldsymbol{\alpha} \in \C^n$, where $0 \leq k \leq n-1$ and $\omega_0=\omega_n=0$.
	\end{enumerate}
	Here, as a vector space $\mathfrak{L}(c, \lambda_1, \lambda_2, \boldsymbol{\alpha}, \psi)$ is equal to $\mathfrak{L}(c, \lambda_1, \lambda_2, \boldsymbol{\alpha})$ with the action defined by $X(b).v=\psi(b)X.v$ for all $X \in \uptau, b \in B, v \in V$ and $\psi: B \rightarrow \C$ is an algebra homomorphism.
\end{theorem}

\begin{proof}
	The result follows directly from Theorem \ref{fg}.
\end{proof}
\section{Classification of Irreducible non-uniformly bounded Harish-Chandra modules over $\uptau(B)$}

In this section, we classify all irreducible Harish-Chandra modules that are not uniformly bounded. First, we define highest weight modules for $\uptau$. Let $\uptau = \oplus_{{\bf{m}} \in \Z^n} \uptau_{\bf{m}}$ be the graded decomposition given by \eqref{hj}. Let $M$ be a subgroup of $\Z^n$ and $\boldsymbol{\beta} \in (\Z^n)^*$ such that $\Z^n= M \oplus \Z \boldsymbol{\beta}$. Define 

\begin{center}
	$\uptau^{-}_M= \bigoplus\limits_{\substack{{\bf{m}} \in M \\ r \in \N}} \uptau_{{\bf{m}}-r\boldsymbol{\beta}}, \quad \uptau_M= \bigoplus\limits_{{\bf{m}} \in M} \uptau_{{\bf{m}}}, \quad \text{and}  \quad \uptau^{+}_M= \bigoplus\limits_{\substack{{\bf{m}} \in M \\ r \in \N}} \uptau_{{\bf{m}}+r\boldsymbol{\beta}}$
\end{center}
Then, $\uptau = \uptau^{-}_M \oplus \uptau_M \oplus \uptau^{+}_M$ is a triangular decomposition of  $\uptau$. Using this decomposition, we define highest weight $\uptau$-modules. Let $X$ be an irreducible $\uptau_M$-module. We extend $X$ to a module over $\uptau_M \oplus \uptau^{+}_M$ by setting $\uptau^{+}_M \cdot X=0$. Then the generalized Verma module for $\uptau$ is defined by 
\begin{center}
	$\mathbb{M}_{\uptau}(X,\boldsymbol{\beta}, M)= \textrm{Ind}_{\uptau_M \oplus \uptau^{+}_M}^{\uptau} X= \mathcal{U}(\uptau) \otimes_{\uptau_M \oplus \uptau^{+}_M}  X.$ 
\end{center}
By standard arguments, $\mathbb{M}_{\uptau}(X,\boldsymbol{\beta}, M)$ has a unique irreducible quotient, denoted by $\mathbb{L}_{\uptau}(X,\boldsymbol{\beta}, M)$. Since $\uptau$ is $\Z^n$-extragraded,  by a result of Billig and Zhao (\cite{bz}, Theorem 1.5), it is known that if $X$ is a uniformly bounded $\Z^{n-1}$-graded exp-polynomial $\uptau_{M}$-module, then 
$\mathbb{L}_{\uptau}(X,\boldsymbol{\beta}, M)$ becomes a Harish-Chandra module for $\uptau$ (see \cite{bz} for definitions).
\begin{theorem}[\cite{sou}, Theorem 6.1]\label{m}
	Any non-trivial irreducible GHW $\uptau$-module $V$ is isomorphic to $\mathbb{L}_{\uptau}(X,\boldsymbol{\beta}, M)$ for some $X$,\,$\boldsymbol{\beta}$,\,and $M$ defined as above.
\end{theorem}

Analogously, we define highest weight modules for $\uptau(B)$. Let $\uptau(B)= \oplus_{{\bf{m}} \in \Z^n} \uptau(B)_{\bf{m}}$ be the graded decomposition given by \eqref{ed}. As before, let $M$ be a subgroup of $\Z^n$ and $\boldsymbol{\beta} \in (\Z^n)^*$ such that $\Z^n= M \oplus \Z \boldsymbol{\beta}$. Similarly define

\begin{center}
	$\uptau(B)^{-}_M= \bigoplus\limits_{\substack{{\bf{m}} \in M \\ r \in \N}} \uptau(B)_{{\bf{m}}-r\boldsymbol{\beta}}, \quad \uptau(B)_M= \bigoplus\limits_{{\bf{m}} \in M} \uptau(B)_{{\bf{m}}}, \quad \text{and}  \quad \uptau(B)^{+}_M= \bigoplus\limits_{\substack{{\bf{m}} \in M \\ r \in \N}} \uptau(B)_{{\bf{m}}+r\boldsymbol{\beta}}$
\end{center}
This gives a triangular decomposition $\uptau(B)= \uptau(B)^{-}_M \oplus \uptau(B)_M \oplus \uptau(B)^{+}_M$. Let $X$ be an irreducible $\uptau(B)_M$-module. Extend $X$ to a module for $\uptau(B)_M \oplus \uptau(B)^{+}_M$ by setting $\uptau(B)^{+}_M \cdot X=0$. Then, the generalized Verma module for $\uptau(B)$ is defined as 
\begin{center}
	$\mathbb{M}_{\uptau(B)}(X,\boldsymbol{\beta}, M)= \textrm{Ind}_{\uptau(B)_M \oplus \uptau(B)^{+}_M}^{\uptau(B)} X= \mathcal{U}(\uptau(B)) \otimes_{\uptau(B)_M \oplus \uptau(B)^{+}_M}  X.$ 
\end{center}
  By standard arguments, $\mathbb{M}_{\uptau(B)}(X,\boldsymbol{\beta}, M)$ has a unique irreducible quotient denoted by $\mathbb{L}_{\uptau(B)}(X,\boldsymbol{\beta}, M)$. The following theorem provides a classification of irreducible non-uniformly bounded Harish-Chandra modules over $\uptau(B)$.
  
  \begin{theorem}\label{oi}
  	Let $V$ be an irreducible Harish-Chandra module for $\uptau(B)$ with non uniformly bounded weight spaces. Then $V$ is a highest weight module for $\uptau(B)$. In particular, 
  	\begin{center}
  		$V \cong \mathbb{L}_{\uptau(B)}(X,\boldsymbol{\beta}, M)$,
  	\end{center} for some $X,\boldsymbol{\beta}$, and  $M$ as above.
  \end{theorem}
\begin{proof}
	In order to prove $V$ is a highest weight module for $\uptau(B)$, we need to find a highest weight vector $v$ such that $\uptau(B)^+v=0$. 
	The proof is exactly analogous to that of \cite[Theorem 4.7]{SPRE}. For completeness and clarity, we summarize the key steps here. One can look at \cite[Theorem 4.7]{SPRE} for the detailed justifications of the steps given below.
	\par
	From \cite[Theorem 4.3]{sp}, there exists a GHW vector $v \in V_{\mu}$ for $\uptau$. Let $U:=\mathcal{U}(\uptau)v$, and define $\mathring{U}=\{u\in U: \uptau.u=0 \}$. Then the quotient $\overline{U}=U/\mathring{U}$ contains no non-zero trivial vectors. By \cite[Lemma 3.6]{maz}, $\overline{U}$ is not uniformly bounded, and hence by Theorem~\ref{m}, it contains an irreducible subquotient of the form $\mathbb{L}_{\uptau}(X,\boldsymbol{\beta}, M)$. Let $\bar{v}$ denote the image of $v$ in $\overline{U}$. Since every element of  $\overline{U}$ is a generalized highest weight vector, by replacing $\bar{v}$ with a suitable vector if necessary, we may assume that
	\begin{center}
		$\mu - \mathbb{Z}_+ \boldsymbol{\beta} + M \setminus \{ \mathbf{0} \} \subseteq \mathrm{Supp}\,(\overline{U})$.
	\end{center}
	Then it follows that $\left\{ \mu + l\boldsymbol{\beta} + M : l \in \mathbb{N} \right\} \cap \mathrm{Supp}(\overline{U}) \neq \emptyset$ for only finitely many values of $l$. Let $L$ be the largest integer such that $\left\{ \mu + L\boldsymbol{\beta} + M \right\} \cap \, \mathrm{Supp}(\overline{U}) \neq \emptyset$. Then 
	\begin{equation*}
		\overline{U}[\mu + L\boldsymbol{\beta} + M]:= \bigoplus_{{\bf{m}} \in M}\overline{U}_{\mu + L\boldsymbol{\beta} + {\bf{m}}}
	\end{equation*}
	 is a uniformly bounded $\uptau_M$-module. Hence from the Lemma 4.2 of \cite{sou}, $\overline{U}[\mu + L\boldsymbol{\beta} + M]$ has a non-zero irreducible $\uptau_M$-submodule $Y$. Thus, $\overline{U}$ contains a highest weight $\uptau$-module given by $P=\mathcal{U}(\uptau)Y$.\par
 Now we need to deal with two cases: (i) $\mathbf{0} \in \mathrm{Supp}(\overline{U})$, (ii) $\mathbf{0} \notin \mathrm{Supp}(\overline{U})$. In both cases, by following the argument of the proof of Theorem 4.7 of \cite{SPRE} line by line, we obtain a highest weight vector for $\uptau(B)$ and that completes the proof.  
	\end{proof}
	We know that if $\mathbb{L}_{\uptau(B)}(X,\boldsymbol{\beta}, M)$ is an irreducible Harish-Chandra module, then $X$ has to be an irreducible cuspidal module for $\uptau(B)_M$ (\cite{bz}, Theorem 1.5). Our next aim is to show that $X$ is a single point evaluation module for $\uptau(B)_M$.

\begin{prop}\label{sd}
	$X$ is a single point evaluation module for $\uptau(B)_M$.
\end{prop}
\begin{proof} Let $\mathring{\uptau}$ be the subalgebra of $\uptau$ generated by the elements 
	\begin{center}
	$\{x \otimes t^{\lm},\, t^{\lm}K_i,\, t^{\lm}d_i: {\lm} \in \Z^{n-1},\,\,  x \in \gg,\, i=1,2,\ldots,n\} $.
	\end{center} 
	Let $\{{\bm{\al}}_2,\bm{\al}_3,\ldots,\bm{\al}_n\}$ be a $\Z$-basis for $M$, and $\bm{\al}_1 \in \Z^n$ such that $\Z^n= \Z\bm{\al}_1 \oplus M $. Let $\{\bm{\beta}_1,\bm{\beta}_2,\ldots, \bm{\beta}_n\}$ be the basis of $\mathbb{R}^n$ that is dual to $\{\bm{\al}_1,\bm{\al}_2,\ldots,\bm{\al}_n\}$. Define $\Phi: \mathring{\uptau} \rightarrow \uptau_M$ by
	\begin{center}
		$x \otimes t^{\lm} \mapsto x \otimes t^{\sum_{i=2}^{n}m_i\bm{\al}_i},\hspace{0.4cm} t^{\lm}K_i \mapsto K(\bm{\al}_i, \sum_{i=2}^{n}m_i\bm{\al}_i) , \hspace{0.4cm} t^{\lm}d_i \mapsto D(\bm{\beta}_i, \sum_{i=2}^{n}m_i\bm{\al}_i) $
	\end{center}
One can easily verify that $\Phi$ is a Lie algebra isomorphism. Thus, we may view $X$ as an irreducible uniformly bounded module for $\mathring{\uptau}(B)$. Actually, we use the subalgebra $\mathring{\uptau}(B)$ instead of $\uptau(B)_M$ for simplicity of calculation. Then, using the same techniques as in Propositions \ref{va},\ref{vb},\ref{vc}, we conclude that $X$ is a single point evaluation module for $\mathring{\uptau}(B)$, and hence for $\uptau(B)_M$.
	\end{proof}

\begin{theorem}\label{sx}
	$\mathbb{L}_{\uptau(B)}(X,\boldsymbol{\beta}, M)$ is a single point evaluation module for $\uptau(B)$.
\end{theorem}
 \begin{proof}
 	From the previous proposition, we know that $X$ is a single point evaluation module for $\uptau(B)_M$. Therefore, there exists an algebra homomorphism $\psi: B \rightarrow \C$ such that
 	\begin{center}
 		 $D({\bf{p}}, {\bf{m}})(b).v=\psi(b)D({\bf{p}}, {\bf{m}}).v$, \hspace{0.5cm} $K({\bf{p}}, {\bf{m}})(b).v=\psi(b)D({\bf{p}}, {\bf{m}}).v$\\
 		 \vspace{0.2cm}
 		 $x \otimes t^{{\bf{m}}}(b)v=\psi(b) x \otimes t^{{\bf{m}}}.v$
 	\end{center}
  for all ${\bf{p}} \in \C^n, {\bf{m}} \in M, b \in B$, and $v \in X$. Let ${\bf{r}} \in \Z^n$ such that ${\bf{r}}=j\boldsymbol{\beta} +  {\bf{m}}$ with $j \in \Z$ and ${\bf{m}} \in M$. By the construction of $\mathbb{L}_{\uptau(B)}(X,\boldsymbol{\beta}, M)$, we know that $\uptau(B)^+$ acts trivially on $X$. Hence, for
 all $j\geq 0,\, {\bf{p}} \in \C^n,\, {\bf{m}} \in M,\, b \in B$, and $v \in X$, we have
 \begin{equation}\label{h}
 D({\bf{p}}, j\boldsymbol{\beta} +  {\bf{m}})(b-\psi(b)).v=K({\bf{p}}, j\boldsymbol{\beta} +  {\bf{m}})(b-\psi(b)).v=x \otimes t^{j\boldsymbol{\beta} + {\bf{m}}}(b-\psi(b))v=0.
 \end{equation}
 Now consider the set, $W=\{v \in V: Y(b).v=\psi(b)\, Y.v=0, \,\,\,  \forall\,\, Y \in \uptau,\, b \in B\}$.\vspace{0.15cm}\\
 {\bf{Claim.}} $X \subset W$. \vspace{0.1cm}\\
 To establish the above claim, we need to prove that \eqref{h} also holds for $j<0$. We use mathematical induction on $j$. For the base case $j=-1$, let $
 {\bf{r}}'=j'\boldsymbol{\beta} +  {\bf{m}}' $ with $j' \in \N$ and ${\bf{m}}' \in M$, then
\begin{align*}
 	\Big(x \otimes t^{{\bf{r}}'}(b')\Big) D({\bf{p}},{\bf{r}})(b-\psi(b)).v &=-({\bf{p}}|{\bf{r}}')x \otimes t^{{\bf{r}+\bf{r}}'}(b'(b-\psi(b))).v\\
 	\Big(D({\bf{p}'}, {\bf{r}'})(b')\Big) D({\bf{p}},{\bf{r}})(b-\psi(b)).v&=\Big(D({\bf{w}}, {\bf{r}+\bf{r}}')
 	+cK({\bf{r}'}, {\bf{r}+\bf{r}}')
 	\Big)(b'(b-\psi(b))).v\\
 	\Big(K({\bf{p}'}, {\bf{r}'})(b')\Big)D({\bf{p}}, {\bf{r}})(b-\psi(b)).v&=K({\bf{w'}}, {\bf{r}+\bf{r}}'))(b'(b-\psi(b))).v,
 \end{align*}
where ${\bf{w}}=({\bf{p'}}|{\bf{r}}){\bf{p}}-({\bf{r'}}|{\bf{p}}){\bf{p'}},\,\, {\bf{w'}}=({\bf{p'}}|{\bf{r}}){\bf{p}}+({\bf{p'}}|{\bf{p}}){\bf{r'}}$ and $c$ is a scalar that depends on $\phi$. Now, if $j'=1$, then ${\bf{r}+\bf{r}}' \in M$, so the above expressions will become zero as $\psi(b'(b-\psi(b)))=\psi(b'b)-\psi(b'\psi(b))=0$. If $j'>1$, then ${\bf{r}+\bf{r}}' \in \uptau(B)_M^+$ and again the above expressions will become zero. Thus $D({\bf{p}},{\bf{r}})(b-\psi(b)).v$ is a highest weight vector which does not lie in the highest weight space, hence must be zero. Similarly, we can show that
\begin{center} 
$K({\bf{p}}, j\boldsymbol{\beta} +  {\bf{m}})(b-\psi(b)).v=0=x \otimes t^{j\boldsymbol{\beta} + {\bf{m}}}(b-\psi(b))v$, \hspace{0.5cm}with \hspace{0.2cm} $j=-1$.
\end{center}
It is now straightforward to verify by induction that the above identity \eqref{h} holds for all $j \leq -1$. It is now easy to check that $W$ is a $\uptau(B)$-submodule of $\mathbb{L}_{\uptau(B)}(X,\boldsymbol{\beta}, M)$. Hence $W=\mathbb{L}_{\uptau(B)}(X,\boldsymbol{\beta}, M)$.
 \end{proof} 
As a consequence of the above theorem, we have $\mathbb{L}_{\uptau(B)}(X,\boldsymbol{\beta}, M) \cong \mathbb{L}_{\uptau}(X,\boldsymbol{\beta}, M)$ with the $\uptau(B)$-action: 
\begin{center}
	$X(b)v=\psi(b)Xv$ \quad for all \,$X \in \uptau,\, b\in B,\, v \in \mathbb{L}_{\uptau(B)}(X,\boldsymbol{\beta}, M)$.
\end{center}
We denote this module by $\mathbb{L}_{\uptau}(X,\boldsymbol{\beta}, M, \psi)$. Now we concentrate on the highest weight space $X$. By Proposition \ref{sd}, we have $\uptau_M \cong \mathring{\uptau}$, so we can consider $X$ as a module over $\mathring{\uptau}$. The next proposition gives us the complete classification of $X$.
\begin{theorem} (\cite{sou}, Theorem 5.2)\label{xs}
	Let $X$ be an irreducible cuspidal $\mathring{\uptau}$-module. Then we have the following.
	\begin{enumerate}
		\item Suppose both $K_1$ and $d_1$ act trivially on $X$. Then for all $\mathbf{m} \in \mathbb{Z}^n$, the operators $t^{\mathbf{m}} K_1$ and $t^{\mathbf{m}} d_1$ also act trivially on $X$. Consequently, $X$ is an irreducible cuspidal module over the full toroidal Lie algebra in the variables $t_2, \ldots, t_n$, which has already been classified in Theorem \ref{fg}.
		
		\item If either $K_1$ or $d_1$ acts non-trivially on $X$, then there exists a finite-dimensional irreducible $\mathfrak{g}$-module $V_1$, a finite-dimensional irreducible $\mathfrak{gl}_{n-1}$-module $V_2$, ${\underline{\boldsymbol{\alpha}}}= (\alpha_2, \ldots, \alpha_n) \in \mathbb{C}^{n-1}$ and $(a, b) \in \mathbb{C}^2 \setminus \{(0, 0)\}$ such that
		\begin{center}
			$	X \cong V_1 \otimes V_2 \otimes \mathbb{C}[t_2^{\pm1}, \ldots, t_n^{\pm1}]$
		\end{center}
		where the action of $\mathring{\uptau}$ on $X$ is defined by:
		
		\begin{center}
			$t^{\lm} K_1(v_1 \otimes v_2 \otimes t^{\lr}) = a (v_1 \otimes v_2 \otimes t^{\lr + \lm})$, \hspace{2.5cm} 
			$t^{\lm} K_i (v_1 \otimes v_2 \otimes t^{\lr}) = 0$, \\
			\vspace{0.45cm}
			$x \otimes t^{\lm}(v_1 \otimes v_2 \otimes t^{\lr})\hspace{-0.1cm}=\hspace{-0.1cm}\,xv_1 \otimes v_2 \otimes t^{\lr + \lm}$,\,\,\,\,\,\, $t^{\lm}\, d_1 (v_1 \otimes v_2 \otimes t^{\lr})\hspace{-0.1cm} = b (v_1 \otimes v_2 \otimes t^{\lr + \lm})$,\\
				\vspace{0.45cm}
			$t^{\lm}\, d_i (v_1 \otimes v_2 \otimes t^{\lr}) = (\alpha_i + r_i)(v_1 \otimes v_2 \otimes t^{\lr + \lm}) 
			+ \sum_{j=2}^n m_j (v_1 \otimes (E_{j,i} v_2) \otimes t^{\lr + \lm})$.
			\end{center}
		\vspace{0.2cm}
		Here, $x \in \mathfrak{g}, \, v_1 \in V_1, \, v_2 \in V_2 , \, \lm,\, \lr \in \mathbb{Z}^{n-1}, \, 2 \leq i \leq n$, and $E_{j,i}$ denotes the matrix of order $n-1$ having 1 at the $(j,i)$-th entry and $0$ elsewhere.
	\end{enumerate}
	
\end{theorem}
\medskip
\begin{remark}\label{yt}
	Observe that the $\mathring{\uptau}$-module $X \cong V_1 \otimes V_2 \otimes \mathbb{C}[t_2^{\pm1}, \ldots, t_n^{\pm1}]$ is completely determined by the parameters $(a,b,c,\lambda_1, \lambda_2, {\underline{\boldsymbol{\alpha}}}) \in \C \times \C \times \C \times P^+_{\gg} \times P^+_{\mathfrak{sl}_n} \times \C^{n-1}$. We denote this module by $\mathfrak{L}(a,b,c, \lambda_1, \lambda_2, {\underline{\boldsymbol{\alpha}}})$. If $(a,b) \neq (0,0)$, $\mathfrak{L}(a,b,c, \lambda_1, \lambda_2, {\underline{\boldsymbol{\alpha}}})$ is an irreducible module over $\mathring{\uptau}$. We continue to use the notation $\mathfrak{L}(a,b,c, \lambda_1, \lambda_2, {\underline{\boldsymbol{\alpha}}})$ even when $(a,b)=(0,0)$, where the $\mathring{\uptau}$-module action remains the same as described in the previous proposition. However, in the case 
	$(a,b)=(0,0)$, the module $\mathfrak{L}(a,b,c, \lambda_1, \lambda_2, {\underline{\boldsymbol{\alpha}}})$ becomes $\mathfrak{L}(c, \lambda_1, \lambda_2, {\underline{\boldsymbol{\alpha}}})$ which can be viewed as a module over the full toroidal Lie algebra $\uptau_{n-1}$ in the variables $t_2,t_3, \ldots, t_n$. This module may not be irreducible over $\uptau_{n-1}$, but it always has a unique irreducible quotient. Hence, for every tuple $(a,b,c,\lambda_1, \lambda_2, {\underline{\boldsymbol{\alpha}}}) \in \C \times \C \times \C \times P^+_{\gg} \times P^+_{\mathfrak{sl}_n} \times \C^{n-1}$, we denote the unique irreducible quotient of $\mathfrak{L}(a,b,c, \lambda_1, \lambda_2, {\underline{\boldsymbol{\alpha}}})$ by $\mathfrak{L}'(a,b,c, \lambda_1, \lambda_2, {\underline{\boldsymbol{\alpha}}})$. 
\end{remark}

\section{Main Theorem}
In this section, we summarize the results obtained so far to classify the irreducible Harish-Chandra modules over the map full toroidal algebra.
\begin{theorem}\label{fc}
	Let $V$ be a non-trivial irreducible Harish-Chandra module over $\uptau(B)$. Then
	\begin{enumerate}
		\item $V$ is either a cuspidal module or a highest weight module. In both cases, $V$ turns out to be a single point evaluation module.
		\item If $V$ is cuspidal, then $V$ is isomorphic to one of the following:
		\begin{enumerate}
			\item $\mathfrak{L}(c, \lambda_1, \lambda_2, \alpha, \psi)$, for some quadruplet $(c, \lambda_1, \lambda_2, \alpha) \in \C \times (P^+_{\gg})^{\times} \times (P^+_{\mathfrak{sl}_n}) \times \C^n$ and $\psi:B \rightarrow \C$.
			\item $\mathfrak{L}(c, 0, \lambda_2, \alpha, \psi)$ for some $(c, \lambda_2, \alpha) \in \C \times (P^+_{\mathfrak{sl}_n}) \times \C^n$ with $(c, \lambda_2) \neq (k, \omega_k)$ for $1 \leq k \leq n-1$ and $\psi:B \rightarrow \C$.
			\item $\mathfrak{L}(n, 0, 0, \alpha, \psi)$ with $\alpha \notin \Z^n$ and $\psi:B \rightarrow \C$.
			\item A submodule $d_K\mathfrak{L}(k, 0, \omega_k, \alpha, \psi)$ of $\mathfrak{L}(k+1, 0, \omega_{k+1}, \alpha, \psi)$ for some $\alpha \in \C^n$, where $0 \leq k \leq n-1$ and $\omega_0=\omega_n=0$ and $\psi:B \rightarrow \C$.
		\end{enumerate}
		\item If $V$ is of highest weight type, then $V$ is isomorphic to $\mathbb{L}_{\uptau}(X,\boldsymbol{\beta}, M, \psi)$ for some $X,\boldsymbol{\beta},M$ as defined in Section 4 and $\psi: B \rightarrow \C$. Moreover, $X$ is isomorphic to $\mathfrak{L}'(a,b,c, \lambda_1, \lambda_2, \underline{\boldsymbol{\alpha}})$ for some $(a,b,c,\lambda_1, \lambda_2, \underline{\boldsymbol{\alpha}}) \in \C \times \C \times \C \times P^+_{\gg} \times P^+_{\mathfrak{sl}_n} \times \C^{n-1}$.
	\end{enumerate}
	\begin{proof}
		\begin{enumerate}
			\item Follows from \Cref{oi}.
			\item Follows directly from \Cref{ds}.
			\item Follows from \Cref{oi},\,\Cref{sx}, \Cref{xs}, and \Cref{yt}.
		\end{enumerate}
	\end{proof}
\end{theorem}
\begin{remark}
	We have come across \cite{bba} after we completed the paper. In that work, the authors study representations of a similar Lie algebra under the additional assumption of integrability. In contrast, we do not assume integrability, which makes our results stronger and more general. Furthermore, they consider weight spaces with respect to the full Cartan subalgebra, whereas we work with weight spaces with respect to the degree-zero derivations. Additionally, the methods employed in our paper are different from those used in \cite{bba}.
\end{remark}
\subsection*{Funding:} The author gratefully acknowledges support from the National Board for Higher Mathematics through a post-doctoral fellowship (Ref. No. 0204/27/(34)/2023/R\&D-II/11935.


\begin{thebibliography}{999999}

		\bibitem[B06]{yb} Billig, Yuly (2006): A category of modules for the full toroidal {L}ie algebra, {\it Int. Math. Res. Not.} Art. ID 68395, 46.
	
	\bibitem[B99]{yb2} Billig Yuly (1999): An extension of the KdV hierarchy arising from a representation of a toroidal Lie algebra, {\it J. Algebra} 217: 40-64.

	\bibitem[BB]{by} Berman, Stephen and Billig, Yuly (1999): Irreducible representations for toroidal {L}ie algebras, {\it J. Algebra} 221(1): 188--231.
	
	\bibitem[BBa]{bba} Pradeep Bisht, Punita Batra (2024): Integrable Modules of Map full Toroidal Lie Algebras, preprint, version 1, arXiv:2410.04495v1.
\bibitem[BF]{bill} Yuly Billig and Vyacheslav Futorny (2016): Classification of irreducible representations of Lie algebra of vector fields on a torus, {\it J. Reine Angew. Math.}\,720: 199–216.
	
	\bibitem[BZ]{bz} Billig, Yuly and Zhao, Kaiming (2004): Weight modules over exp-polynomial {L}ie algebras {\it J. Pure Appl. Algebra} 191(1-2): 23--42.
	

	
	\bibitem[CLW]{CLW} Cai, Yan-an and L\"u, Rencai and Wang, Yan (2021): Classification of simple {H}arish-{C}handra modules for map
	(super)algebras related to the {V}irasoro algebra, {\it J. Algebra }, 570: 397--415.
	
	\bibitem[F]{F} Futorny, V. M. (1996): Irreducible non-dense {$A^{(1)}_1$}-modules. {\it Pacific J. Math.} 172(1): 83--99.
	
	
	\bibitem[GHL]{gao} Gao, Yun and Hu, Naihong and Liu, Dong (2016): Representations of the affine-{V}irasoro algebra of type {$A_1$}, {\it J. Geom. Phys.}, 106: 102--107.
	
	\bibitem[GZ]{guo11} Xiangqian Guo and Kaiming Zhao (2011): \newblock Irreducible weight modules over Witt algebras, {\it Proc. Amer. Math. Soc.}\,139: 2367–2373.
	
	\bibitem[GLZ] {ZRX} Guo, Xiangqian and Lu, Rencai and Zhao, Kaiming (2011): Simple {H}arish-{C}handra modules, intermediate series-modules, and {V}erma modules over the loop-{V}irasoro algebra. {\it Forum Math.} 23(5): 1029--1052.
	
	\bibitem[GLZ14]{guo} Xiangqian Guo, Gengqiang Liu, and Kaiming Zhao (2014): Irreducible Harish-Chandra modules over extended Witt algebras, {\it Ark. Mat.} \,52(1): 99–112.
	


	\bibitem[JY]{ji} Jiang, Cuibo and You, Hong (2004) Irreducible representations for the affine-{V}irasoro {L}ie algebra of type {$B_l$}, {\it Chinese Ann. Math. Ser. B}, 25(3): 359--368.
	
		\bibitem[K]{kl} Kassel, C. (1984) : K\"{a}hler differentials and coverings of complex simple {L}ie algebras extended over a commutative algebra, {\it J. Pure Appl. Algebra} 34(2-3):265--275.
	
	\bibitem[Kac]{ka} Kac, Victor G (1986): Highest weight representations of conformal current algebras. In Topological and geometrical methods in field theory ({E}spoo, 1986). World Sci. Publ., Teaneck, NJ, pp. 3--15.
	
	\bibitem[KR]{bom} Kac, V.G., Raina, A.K. (1988): Highest weights representations of infinite dimensional Lie algebras. Adv. Ser. Math. Phys. 2.
	
		\bibitem[L]{ta} Larsson, T. A.(1992): Conformal fields: a class of representations of {${\rm Vect}(N)$} {\it Internat. J. Modern Phys. A}, 26: 6493--6508.
	
	\bibitem[LQ]{bo} Liu, Xu Feng and Qian, Min (1994): Bosonic {F}ock representations of the affine-{V}irasoro algebra, {\it J. Phys. A} 27(5): L131--L136G.
	
	\bibitem[MRY]{yokunama} Moody, Robert V,  Eswara Rao, S and Yokonuma, Takeo (1990): Toroidal {L}ie algebras and vertex representations. {\it Geom. Dedicata} 35(1-3):283--307.
	
	\bibitem[MZ]{maz} Mazorchuk, V and Zhao, K (2011): Supports of weight modules over {W}itt algebras. {\it Proc. Roy. Soc. Edinburgh Sect. A} 141(1): 155--170.
	
	\bibitem[M]{om} Mathieu, Olivier (1992): Classification of {H}arish-{C}handra modules over the {V}irasoro {L}ie algebra, {\it Invent. Math.} 107(2):225--234.
	
	\bibitem[P]{sou} Pal Souvik (2022): Classification of irreducible Harish-Chandra modules over full toroidal Lie algebras and higher-dimensional Virasoro algebras, preprint, version 3, arXiv:2203.06148v3.
	
	\bibitem[PR]{sp} Pal, Souvik and Eswara Rao, S (2021): Classification of level zero irreducible integrable modules for twisted full toroidal {L}ie algebras {\it J. Algebra} 578: 1--29.
	
	\bibitem[R96]{rao2} S.~Eswara Rao (1996): \newblock Irreducible representations of the Lie-algebra of the diffeomorphisms of a $d$-dimensional torus, {\it J. Algebra}\,182: 401–421. 
	
	\bibitem[R04]{2004} Eswara Rao, S (2004): Classification of irreducible integrable modules for toroidal {L}ie algebras with finite dimensional weight spaces. {\it J. Algebra} 277(1):318-348.
	
	\bibitem[RSM]{fi} Eswara Rao, S and Sharma, Sachin S. and Mukherjee, S (2021): Integrable modules for loop affine-{V}irasoro algebras, {\it Comm. Algebra} 49(12): 5500--5512.
	
		\bibitem[S]{SAV} Savage, A (2012): Classification of irreducible quasifinite modules over map {V}irasoro algebras. {\it Transform. Groups} 17(2): 547--570.
	
	
	\bibitem[SCPR]{SPRE} Sharma, Sachin S. and Chakraborty, Priyanshu and Pandey, Ritesh Kumar and Eswara Rao, S.(2024): Representations of map extended {W}itt algebras {\it J. Algebra} 639: 327--353.
	
		\bibitem[Sh]{sh} Shen, Guang Yu(1986): Graded modules of graded {L}ie algebras of {C}artan type. {I}. {M}ixed products of modules, {\it Sci. Sinica Ser. A}, 29(6): 570--581.
	
	
	

	
	
	

		
	
	
	
	
	\end{thebibliography}
\end{document}